\documentclass[12pt, 14paper,reqno]{amsart}
\usepackage{amsmath,amsfonts,amssymb}
\usepackage{mathrsfs}
\usepackage{longtable}
\usepackage[breaklinks]{hyperref}
\usepackage{graphicx}
\makeatletter
\@namedef{subjclassname@2020}{%
  \textup{2020} Mathematics Subject Classification}
\makeatother

\theoremstyle{plain}
\newtheorem{thm}{Theorem}[section]
\newtheorem{lem}{Lemma}[section]
\newtheorem{cor}{Corollary}[section]

\newtheorem{conj}{Conjecture}[section]

\newtheorem{remark}{Remark}[section]

\newtheorem{thma}{Theorem}

\theoremstyle{proof}

\numberwithin{equation}{section}

\begin{document} 
\title[On a conjecture of Iizuka]{On a conjecture of Iizuka}
\author{Azizul Hoque}

\email{ahoque.ms@gmail.com}

\subjclass[2020]{11R11; 11R29}

\date{\today}

\keywords{Imaginary quadratic field, Class number, Iizuka's conjecture, Exponent}

\begin{abstract}
For a given odd positive integer $n$ and an odd prime $p$, we construct an infinite family of quadruples of imaginary quadratic fields $\mathbb{Q}(\sqrt{d})$, $\mathbb{Q}(\sqrt{d+1})$, $\mathbb{Q}(\sqrt{d+4})$ and $\mathbb{Q}(\sqrt{d+4p^2})$ with $d\in \mathbb{Z}$ such that the class number of each of them is divisible by $n$. Subsequently, we show that there is an infinite family of quintuples of imaginary quadratic fields $\mathbb{Q}(\sqrt{d})$, $\mathbb{Q}(\sqrt{d+1})$, $\mathbb{Q}(\sqrt{d+4})$, $\mathbb{Q}(\sqrt{d+36})$ and $\mathbb{Q}(\sqrt{d+100})$ with $d\in \mathbb{Z}$ whose class numbers are all divisible by $n$. Our results  provide a complete proof of Iizuka's conjecture (in fact a generalization of it) for the case $m=1$.  Our results also affirmatively answer a weaker version of (a generalization of) Iizuka's conjecture for $m\geq 4$. 

\end{abstract}
\maketitle{}

\section{Introduction}
It has been proved that there are infinitely many real (resp. imaginary) quadratic fields with class numbers divisible by a given positive integer (see \cite{AC55, CH18, SO00, YA70}). An analogous problem for tuples of quadratic fields arises from Scholz's Spiegelungssatz \cite{SC32}. In \cite{KO02}, Komatsu  studied this problem for a pair of quadratic fields and proved that there are infinitely many pairs of quadratic fields $\mathbb{Q}(\sqrt{d})$ and $\mathbb{Q}(\sqrt{md})$  with $m, d\in \mathbb{Z}$ whose class numbers are divisible by $3$. Later, he generalized this result in \cite{KO17} to $n$-divisibility of the class numbers of pairs of imaginary quadratic fields. On the other hand, Iizuka \cite{IZ18} studied a slight variant of this problem and construct an infinite family of pairs of imaginary quadratic fields $\mathbb{Q}(\sqrt{d})$  and $\mathbb{Q}(\sqrt{d+1})$ with $d\in\mathbb{Z}$ whose class numbers are divisible by $3$. Further, he posed the following conjecture in the same paper.

\begin{conj}\label{conj1}
For any prime number $p$ and any positive integer $m$, there is an infinite family of $m + 1$ successive real (or imaginary) quadratic fields 
$$\mathbb{Q}(\sqrt{d}), \mathbb{Q}(\sqrt{d+1}), \cdots, \mathbb{Q}(\sqrt{d+m})$$
with $d\in \mathbb{Z}$ whose class numbers are divisible by $p$.
\end{conj}

In \cite{CM21}, Chattopadhyay and Muthukrishnan extended Iizuka's result from pairs to certain triples of imaginary quadratic fields following the methods used in \cite{IZ18}. In other words, they gave an affirmative answer of a weaker version of Conjecture \ref{conj1} for $p = 3$. It follows from a recent result of Iizuka, Konomi and Nakano  \cite{IK21} that Conjecture \ref{conj1} is true for $m=1$ when $p=3,5,7$. 

Very recently, Krishnamoorthy and Pasupulati \cite{KP21} cleverly used \cite[Theorem 1]{LO09} and an extended version of \cite[Theorem 3.2]{HC19} to settled Conjecture \ref{conj1} for $m=1$. 

In this paper, for a given odd prime $p$ we construct an infinite family of quadruples of imaginary quadratic fields $\mathbb{Q}(\sqrt{d})$, $\mathbb{Q}(\sqrt{d+1})$, $\mathbb{Q}(\sqrt{d+4}) $ and $\mathbb{Q}(\sqrt{d+4p^2})$ with $d\in\mathbb{Z}$ whose class numbers are all divisible by a given odd integer $n\geq 3$. This extends the results of \cite{CM21, IK21, KP21} in both directions; from pairs to quadruples/quituples of fields and from primes to odd integers. It also gives a proof of a weaker version of Conjecture \ref{conj1} for any prime $p\geq 3$ (in fact for any odd integer $n\geq 3$). The precise statement of our first result is the following:
\begin{thm}\label{thm}
For any odd positive integer $n$ and any odd prime $p$, there are infinitely many quadruples of imaginary quadratic fields $\mathbb{Q}(\sqrt{d})$, $\mathbb{Q}(\sqrt{d+1})$, $\mathbb{Q}(\sqrt{d+4})$ and $\mathbb{Q}(\sqrt{d+4p^2})$ whose class numbers are all divisible by $n$.
\end{thm} 
Note that we can construct an infinite family of quintuples or higher tuples of imaginary quadratic fields whose class numbers are all divisible by $n$ by choosing different values for the prime $p$. For instance, utilizing Corollary \ref{corx} and following the proof of Theorem \ref{thm}, we get the following:
\begin{thm}\label{thm1}
For a given odd integer $n\geq 1$, there are infinitely many quintuples of imaginary quadratic fields 
$$\left(\mathbb{Q}(\sqrt{d}), \mathbb{Q}(\sqrt{d+1}), \mathbb{Q}(\sqrt{d+4}), \mathbb{Q}(\sqrt{d+36}), \mathbb{Q}(\sqrt{d+100}) \right)$$
with $d\in\mathbb{Z}$ whose class numbers are all divisible by $n$. 
\end{thm} 
Our method relies on the prominent result of Bilu, Hanrot and Voutier \cite{BH01} concerning the primitive divisors of Lehmer numbers as well as on the solvability of certain Lebesgue-Ramanujan-Nagell type equations considered in \cite{HE99}. This method does not allow us to include an imaginary quadratic field of the form $\mathbb{Q}(\sqrt{d+m})$ in the tuple, when $m$ is non-square integer.

\section{$n$-divisibility of the class-numbers of $\mathbb{Q}\left(\sqrt{1-4U^n}\right)$ and $\mathbb{Q}\left(\sqrt{1-V^n}\right)$}
Here, we recall some results concerning the $n$-divisibility of class numbers of the imaginary quadratic fields   $\mathbb{Q}\left(\sqrt{1-4U^n}\right)$ and $\mathbb{Q}\left(\sqrt{1-V^n}\right)$. These results will be used in the proof of Theorem \ref{thm}. 
\begin{thma}\label{thma}
If $n\geq 3$ is an odd integer, then for any integer $U\geq  2$ the ideal class groups of the imaginary quadratic fields $\mathbb{Q}(\sqrt{1-4U^n})$ contain an element of order $n$.
\end{thma} 
In 1978, Gross and Rohrlich gave the outline of a proof of this  theorem (see \cite[Theorem 5.3 and Remark 2]{GR78}). Their method of proof was based upon the affine points on the Fermat curve $x^p+y^p=1$ over the imaginary quadratic field $\mathbb{Q}(\sqrt{1-4U^n})$. Later, Louboutin \cite{LO09} gave a complete proof of this theorem using number theoretic technique. It follows from Siegel's theorem (see \cite{SI29}) that for each integer $d>1$ there are at most finitely many positive integers $U$ such that $1-4U^n=-dX^2$. This ensures the infinitude of the above family of imaginary quadratic fields.      

The $n$-divisibility of the class numbers of the family of imaginary quadratic fields $\mathbb{Q}(\sqrt{1-V^n})$ was studied by Nagell \cite[Theorem 25]{NA22} for any odd integers $V\geq 3$ and $n\geq 3$. Later, Murty \cite[Theorems 1 and 2]{RM98} proved that the class group of the imaginary quadratic field $\mathbb{Q}(\sqrt{1-V^n})$ has an element of order $n$ when either $V^n-1$ is square-free or its square part is $<V^{n/2}/8$ for any odd integers $V\geq 5$ and $n\geq 3$. However, it follows from the fact $DX^2+1=V^n$ has no integer solution when both $V$ and $n$ are odd, except for $(V, n)=(5, 3), (7,3), (13, 3)$ (see \cite{HJ04}), that the above conditions are no longer required. This fact also confirms that there are infinitely many such imaginary quadratic fields. Finally, this was elucidated by Cohn in \cite[Corollary 1]{CO03} as follows:
\begin{thma}\label{thmb} 
Assume that $n\geq 3$ and $V\geq 3$ are odd integers. Then the class number of the imaginary quadratic field $\mathbb{Q}(\sqrt{1-V^n})$ is divisible by $n$, except for $(n,V)=(5,3)$.
\end{thma}

\section{The divisibility of the class number of $\mathbb{Q}(\sqrt{p^2-\ell^n})$}
Many special cases of the divisibility of the class number of the imaginary quadratic field $\mathbb{Q}(\sqrt{x^2-y^n})$ have been studied with some restrictions on $x, y$ and $n$. One of such restrictions is that $y$ is an odd prime (see \cite{CH18} and the references therein), and hence none of the known results can be used to complete the proof of Theorem \eqref{thm}. Thus, we consider a family of imaginary quadratic fields of the above form where $y$ is not a prime and it will be a useful ingredient in the proof of Theorem \ref{thm}. Here, we mainly prove:
 
\begin{thm}\label{thm2}
Let $\ell>1$ and $n>1$ be odd integers, and $p$ an odd prime such that $\ell\equiv 3\pmod 4, ~\gcd(\ell, p)=1$ and $p^2<\ell^n$. Assume that $-d$ is the square-free part of $p^2-\ell^n$. If $p\not\equiv \pm 1\pmod d$, then the class number of $\mathbb{Q}(\sqrt{p^2-\ell^n})$ is divisible by $n$.  
\end{thm}
Theorem \ref{thm2} extends \cite[Theorem 1.1]{CH18}, where the authors assumed that $\ell$ is an odd prime. This primality condition on $\ell$ restricts us to apply \cite[Theorem 1.1]{CH18} in the proof of Theorem \ref{thm}. Further, following the proof of \cite[Theorem 1.2]{CH18}, we make the following remark.
\begin{remark}\label{rem1}
The family of imaginary quadratic fields discussed in Theorem \ref{thm2} has infinitely many members.
 \end{remark}
 Now for $p=3, 5$, the condition `$p\not\equiv \pm 1\pmod d$' can be removed by applying \ref{eqc6} except for $(\ell,n)=(3,3)$, and thus we have the following straightforward corollary.
 \begin{cor}\label{corx}
 Let $\ell>1$ and $n>1$ be as in Theorem \ref{thm2}. For $x=3,5$ with $\gcd(\ell, x)=1$, the class number of $\mathbb{Q}(\sqrt{x^2-\ell^n})$ is divisible by $n$ except the case $(\ell,n)=(3,3)$. 
 \end{cor}

The proof of Theorem \ref{thm2} relies on the prominent result of Bilu, Hanrot and Voutier \cite{BH01} on existence of primitive divisors of Lehmer numbers. 

\subsection{Lehmer numbers and their primitive divisors}
A pair $(\alpha, \beta)$ of algebraic integers is said to be a Lehmer pair if $(\alpha + \beta)^2$ and $\alpha\beta$ are two non-zero coprime rational integers, and $\alpha/\beta$ is not a root of unity. For a given positive integer $n$, the Lehmer numbers correspond to the pair $(\alpha, \beta)$ are defined as 
$$\mathcal{L}_n(\alpha, \beta)=\begin{cases}
\dfrac{\alpha^n-\beta^n}{\alpha-\beta}, & \text{ if $n$ is odd,} \vspace{1mm}\\
\dfrac{\alpha^n-\beta^n}{\alpha^2-\beta^2}, & \text{ if $n$ is even}.
\end{cases}$$
It is known that all Lehmer numbers are non-zero rational integers.  Two Lehmer pairs $(\alpha_1, \beta_1)$ and $(\alpha_2, \beta_2)$ are said to be equivalent if $\alpha_1/\alpha_2=\beta_1/\beta_2\in \{\pm 1, \pm\sqrt{-1} \}$. A prime divisor $p$ of $\mathcal{L}_n(\alpha, \beta)$ is said to be primitive if $p\nmid(\alpha^2-\beta^2)^2
\mathcal{L}_1(\alpha, \beta) \mathcal{L}_2(\alpha, \beta) \cdots \mathcal{L}_{n-1}(\alpha, \beta)$. 
The following classical result  was proved in \cite[Theorem 1.4]{BH01}. 
\begin{thma}\label{thmBH}
The Lehmer number $\mathcal{L}_n (\alpha, \beta) $ has primitive divisors for any integer $n>30$. 
\end{thma}

Given a Lehmer pair $(\alpha, \beta)$, let $a=(\alpha+\beta)^2$ and $b=(\alpha-\beta)^2$. Then $\alpha=(\sqrt{a}\pm\sqrt{b})/2$ and $\beta=(\sqrt{a}\mp\sqrt{b})/2$. The pair $(a, b)$ is called the parameters corresponding to the Lehmer pair $(\alpha, \beta)$. The following lemma is extracted from \cite[Theorem 1]{VO95}.
\begin{lem}\label{lemVO}
Let $t$ be an odd integer such that $7\leq t\leq 29$. If the Lehmer numbers $\mathcal{L}_t(\alpha, \beta)$ have no primitive divisor, then up to equivalence, the parameters $(a, b)$ of the corresponding pair $(\alpha, \beta)$ are as follows:
\begin{itemize}
\item[(i)] $(a, b)=(1,-7), (1, -19), (3, -5), (5, -7), (13, -3), (14, -22)$, when $t=7$;
\item[(ii)] $(a,b) = (5,-3), (7,-1),(7,-5)$, when $t=9$;
\item[(iii)] $ (a, b)=(1,-7)$, when $t=13$;
\item[(iv)] $(a,b) = (7,-1),(10,-2)$, when $t=15$. 
\end{itemize}
\end{lem}
Let $F_k$ (resp. $L_k$) denote the $k$-th term in the Fibonacci (resp. Lucas) sequence defined by $F_0=0,   F_1= 1$,
and $F_{k+2}=F_k+F_{k+1}$ (resp. $L_0=2,  L_1=1$, and $L_{k+2}=L_k+L_{k+1}$), where $k\geq 0$ is an integer. 
The following lemma is a part of \cite[Theorem 1.3]{BH01}. 
\begin{lem}\label{lemBH}
For $p=3, 5$, let the Lehmer numbers $\mathcal{L}_p(\alpha, \beta)$ have no primitive divisor. Then up to equivalence, the parameters $(a, b)$ of the corresponding pair $(\alpha, \beta)$ are:
\begin{itemize}
\item[(i)] For $p=3, (a, b)=\begin{cases}(1+t, 1-3t) \text{ with } t\ne 1,\\
 (3^k+t, 3^k-3t) \text{ with } t\not\equiv 0\pmod 3, (k,t)\ne(1,1);\\ \end{cases}$ 
 \vspace{2mm}
\item[(ii)] For $p=5, (a, b)=\begin{cases}
(F_{k-2\varepsilon}, F_{k-2\varepsilon}-4F_k)\text{ with } k\geq 3,\\
 (L_{k-2\varepsilon}, L_{k-2\varepsilon}-4L_k)\text{ with } k\ne 1;
 \end{cases}
 $
\end{itemize}
where $t\ne 0$ and $k\geq 0$ are any integers and $\varepsilon=\pm 1$.
\end{lem}
\subsection{Two important lemmas}
Given an integer $D\equiv 0, 1\pmod 4$, assume that $h^*(D)$ is the class number of binary quadratic primitive forms with discriminant $D$. Also for a square-free integer $d$, let $h(d)$ denote the class number of $\mathbb{Q}(\sqrt{d})$. Then we have (cf. \cite[\S16.13; p. 444]{HU87} the following: 
\begin{lem}\label{lemf}
Let $d\equiv 2\pmod 4$ be a square-free positive integer. Then
$h(-d)=
h^*(-4d).$
\end{lem}
The following lemma is a special case of \cite[Theorem 6.2]{HE99} when $(D_1,D_2)=(1,-d)$. 
\begin{lem}\label{lems}
Let $d>3$ and $\ell >1$ be integers such that $\gcd(\ell,2d) =1$. If the equation
$$x^2 +dy^2 =\ell^z,~~ x, y, z\in\mathbb{N}, ~\gcd(x,y)=1$$
has a solution, then all the solutions $(x, y, z)$ of this equation can be expressed as 
$$x+y\sqrt{-d}=\varepsilon(a +\mu b\sqrt{-d})^t, ~~z=st,$$
where $a,b,s,t$ are positive integers satisfying
$$a^2+db^2=\ell^s, ~~ \gcd(a,b)=1\text{ and }s\mid h^*(-4d)$$
and $\varepsilon, \mu \in \{- 1, 1\}$.
\end{lem}
\subsection{Proof of Theorem \ref{thm2}}
Let $d$ be the square-free part of $\ell^n-p^2$. Then $p^2-\ell^n=-dr^2$ for some $r\in \mathbb{N}$, and thus $(x,y,z)=(p, r,n)$ is a positive integer solution of the equation 
$$x^2+dy^2=\ell^z, ~\gcd(x,y)=1.$$
Thus by Lemma \ref{lems}, we get
\begin{equation}\label{eqc1}
p+r\sqrt{-d}=\varepsilon(a+\mu b\sqrt{-d})^t,~~ \varepsilon, \mu \in \{-1, 1\}
\end{equation}
with 
\begin{equation}\label{eqc2}
n=st,~~s, t\in \mathbb{N}, 
\end{equation}
where $a$ and $b$ are positive integers satisfying
\begin{equation}\label{eqc3}
a^2+db^2=\ell^{s},~~ \gcd(a,b)=1 
\end{equation}
and 
\begin{equation}\label{eqc4}
s\mid h^*(-4d).
\end{equation}
Since $\ell\equiv 3\pmod 4$ and $n$ is odd, so that $p^2-\ell^n=-dr^2$ gives $d\equiv 2\pmod 4$ and $r$ is odd. Also both $s$ and $t$ are  odd as $n$ is odd. Further reading \eqref{eqc3} modulo $4$, we get $a^2+2b^2\equiv 3\pmod 4$ as $\ell\equiv 3\pmod 4$ and $s$ is odd, which ensures that both $a$ and $b$ are odd. 

We now equate the real parts from both sides in \eqref{eqc1} to get 
\begin{equation}\label{eqc5}
p=\varepsilon a\sum^{\frac{t-1}{2}}_{j=0} \binom{t}{2j} a^{t-2j-1}(-db^2)^j.
\end{equation}
This implies $a\mid p$ and thus $a=1,p$. If $a=1$, then it becomes
\begin{equation}\label{eqc5*}
\sum^{\frac{t-1}{2}}_{j=0} \binom{t}{2j} (-db^2)^j=p\varepsilon =\pm p.
\end{equation}
Reading \eqref{eqc5*} modulo $d$, we get $p\equiv \pm 1\pmod d$, which contradicts to the assumption. 
Therefore $a=p$ and thus \eqref{eqc5} becomes 
\begin{equation}\label{eqc6}
\sum^{\frac{t-1}{2}}_{j=0} \binom{t}{2j} p^{t-2j-1}(-db^2)^j=\varepsilon=\pm 1.
\end{equation}

As $a=p$, so that  \eqref{eqc1} reduces to 
\begin{equation}\label{eqc7}
p+r\sqrt{-d}=\varepsilon(p+\mu b\sqrt{-d})^t,~~ \varepsilon, \mu \in \{-1, 1\}.
\end{equation}
 We now assume that $\alpha=\mu b\sqrt{-d}+p$ and $\beta=\mu b\sqrt{-d}-p$. Then both $\alpha$ and $\beta $ are algebraic integers. Clearly, $(\alpha+\beta)^2=-4db^2$ and $\alpha\beta=-p^2-db^2=-\ell^s$ (by \eqref{eqc3}) are coprime rational integers. Furthermore, it follows from the following identity 
 $$\frac{4db^2}{\ell^s}=\frac{(\alpha+\beta)^2}{\alpha\beta}=\frac{\alpha}{\beta}+\frac{\beta}{\alpha}+2$$
 that 
 $$\ell^s\left( \frac{\alpha}{\beta}\right)^2+2(\ell^s-2db^2)\frac{\alpha}{\beta}+\ell^s=0.$$
 Since $\ell>1$ and $\gcd(\ell^s, 2(\ell^s-2db^2))=\gcd(\ell^s, 4db^2)=\gcd(p^2+db^2, 4db^2)=1$, so that the last equation shows that $\dfrac{\alpha}{\beta}$ is not a root of unity. Therefore $(\alpha, \beta)$ is a Lehmer pair  with parameters $(-4db^2, 4p^2)$ and thus the corresponding Lehmer number for $t$ is 
 $$\mathcal{L}_t(\alpha, \beta)=\frac{\alpha^t-\beta^t}{\alpha-\beta}$$
 as $t$ is odd. Utilizing \eqref{eqc7}, we get $$|\mathcal{L}_t(\alpha, \beta)|=1.$$  
 This confirms that the Lehmer number $\mathcal{L}_t(\alpha, \beta)$ has no primitive divisor, and hence Theorem \ref{thmBH} and Lemma \ref{lemVO} (utilizing the fact that $(-4db^2, 4p^2)$ is the parameters) ensure that $t\in\{1,3,5\}$. 
 
 In case of $t=5$, we get by Lemma \ref{lemBH} that $-4db^2=F_{ k-2\varepsilon}$ or $-4db^2=L_{ k-2\varepsilon}$. Clearly, none of these is possible. 
 
 Finally for $t=3$, \eqref{eqc6}  implies that $p^2-3db^2=\pm 1$. Reading it modulo $4$, we see that `$+$' sign is not possible, and thus $p^2-3db^2=-1$. This is not possible by reading it modulo $3$. 
 
  Therefore $t=1$, and thus \eqref{eqc2} and  \eqref{eqc4} together imply that $n\mid h^*(-4d)$. Thus, we complete the proof by Lemma \ref{lemf}.

%

\section{Proof of Theorem \ref{thm}}
We first fix an odd integer $n\geq 3$. We now define the set
$$\mathcal{N}_n=\left\{ k\in \mathbb{N}:  n\mid h(1-4k^n) \right\}.$$ 
Then by Theorem \ref{thma} the set $\mathcal{N}_n$ is infinite. 

Now for any $k\in \mathcal{N}_n$, we set $d=4(1-4k^n)^n$. Then $\mathbb{Q}(\sqrt{d})=\mathbb{Q}(\sqrt{1-4k^n})$ as $n$ is odd. Thus, by Theorem \ref{thma} there are infinitely many such $d$ satisfying $n\mid h(d)$. In other words, $\mathcal{F}(\mathcal{N}_n)=\{\mathbb{Q}(\sqrt{1-4k^n}):k\in\mathcal{N}_n\}$ is an infinite set. 

Now we assume that $U=4k^n-1$ with $k\in \mathcal{N}_n$. Then $1-4U^n=1-4(4k^n-1)^n=4(1-4k^n)^n+1=d+1$, and thus by Theorem \ref{thma}, we have $n\mid h(d+1)$.  

Again for $k\in \mathcal{N}_n$, let us assume that $V=4k^n-1$. Then $V\geq 3$ and is odd, and thus by Theorem \ref{thmb}, we get $n \mid h(1-V^n)$. Since $4(1-V^n)=4-4(4k^n-1)^n=4+4(1-4k^n)^n=d+4$ and $\mathbb{Q}(\sqrt{4(1-V^n)})=\mathbb{Q}(\sqrt{(1-V^n)})$, so that $n\mid h(d+4)$. 

Finally for any $k\in \mathcal{N}_n$, let $\ell=4k^n-1$. Then $\ell\equiv 3\pmod 4$ and hence by utilizing Theorem \ref{thm2}, we have $n\mid h(p^2-\ell^n)$ for any odd prime $p$ satisfying $p\not\equiv\pm 1\pmod d$. Here, $d$ is the square-free part of $\ell^n-p^2$.  Now $4(p^2-\ell^n)=4p^2-4\ell^n=4p^2-4(4k^n-1)^n=d+4p^2$, which implies that $\mathbb{Q}(\sqrt{d+4p^2})=\mathbb{Q}(\sqrt{4(p^2-\ell^n)})$, and thus $n\mid h(d+4p^2)$. This completes the proof of Theorem \ref{thm}.  
\section{Concluding remarks}
In \cite{XC20}, Xie and Chao studied Conjecture \ref{conj1} and proved the following result using Yamamoto's \cite{YA70} construction. 
\begin{thma}\label{thmc}
For any odd positive integer $n$ and any positive integer $m$, there are infinitely many pairs of imaginary fields $\mathbb{Q}(\sqrt{d})$ and $\mathbb{Q}(\sqrt{d + m})$ whose class groups have an element of order $n$ respectively.
\end{thma}
Theorem \ref{thmc} can be viewed as a weaker variant of a generalization of Conjecture \ref{conj1}.  For $m=1$, it provides a generalization of the main result of \cite{KP21} though \cite{XC20} appeared before \cite{KP21}.  In other words, Theorem \ref{thmc} gives a complete proof of the following generalization of Conjecture \ref{conj1} for $m=1$ and a proof of a weaker version of the same for $m\geq 2$.   
\begin{conj}\label{conj2} 
For any odd integer $n\geq 3$ and any integer $m\geq 1$, there is an infinite family of $m + 1$ successive imaginary (or real) quadratic fields 
$$\mathbb{Q}(\sqrt{d}), \mathbb{Q}(\sqrt{d+1}), \cdots, \mathbb{Q}(\sqrt{d+m})$$
with $d\in \mathbb{Z}$ whose class numbers are all divisible by $n$.
\end{conj}
Theorem \ref{thm} offers a constructive proof of Conjecture \ref{conj2} for $m=1$. This theorem also offers a proof of a weaker version of Conjecture \ref{conj2} for $m=4$, which has missed the families of imaginary quadratic fields $\mathbb{Q}(\sqrt{d+2})$ and $\mathbb{Q}(\sqrt{d+3})$ from the complete proof. When $m=4p^2$ with $p$ an odd prime, Theorem \ref{thm} presents a proof of a weaker version of Conjecture \ref{conj2}. We complete this paper by the following remark.
\begin{remark}
For a given positive integer $m$, let $p_m$ denote the largest prime less than or equal to $m$ and $\pi(m)$ the prime-counting function.  Then for a given positive odd integer $n$, our construction gives an infinite family of at least $(\pi(m)+2)$-tuples of imaginary quadratic fields, $$\left(\mathbb{Q}(\sqrt{d}), \mathbb{Q}(\sqrt{d+1}), \mathbb{Q}(\sqrt{d+4}), \mathbb{Q}(\sqrt{d+36}), \cdots \mathbb{Q}(\sqrt{d+4p_m^2}) \right)$$
with $d\in \mathbb{Z}$ whose class numbers are all divisible by $n$. 
\end{remark}   
\section*{Acknowledgements}
The author is grateful to Professor K. Chakraborty and Professor Yasuhiro Kishi for their valuable comments on the paper. The author is thankful to Professor Y. Iizuka for providing a copy of   \cite{IK21}. 
The author gratefully acknowledges the anonymous referee for his/her valuable remarks that immensely improved the results as well as the presentation of the paper.
This work was supported by SERB-NPDF (PDF/2017/001958), Govt. of India.


\begin{thebibliography}{99}
\bibitem{AC55} N. C. Ankeny and S. Chowla, {\it On the divisibility of the class number of quadratic fields}, Pacific J. Math. {\bf 5} (1955), 321--324.


\bibitem{BH01} Y. Bilu, G. Hanrot and P. M. Voutier, {\it Existence of primitive divisors of Lucas and Lehmer numbers (with an appendix by M. Mignotte)}, J. Reine Angew. Math. {\bf 539} (2001), 75--122.

\bibitem{CH18} K. Chakraborty, A. Hoque, Y. Kishi and P. P. Pandey, {\it  Divisibility of the class numbers of imaginary quadratic fields}, J. Number Theory {\bf 185} (2018), 339--348.
\bibitem{CM21} J. Chattopadhyay and  S. Muthukrishnan, {\it On the simultaneous $3$-divisibility of class numbers of triples of imaginary quadratic fields}, Acta Arith. {\bf 197} (2021), no. 1,  105--110. 
\bibitem{CO03} J. H. E. Cohn, {\it On the Diophantine equation $x^n = Dy^2 + 1$}, Acta Arith. {\bf 106} (2003), no. 1, 73--83.

\bibitem{GR78} B. H. Gross and D. E. Rohrlich, {\it Some results on the Mordell-Weil group of the Jacobian of the Fermat curve}, Invent. Math. {\bf 44} (1978), 201--224.
\bibitem{HJ04} E. Herrmann, I. J\'ar\'asi and A. Peth\H{o}, {\it Note on J. H. E. Cohn's paper ``The Diophantine equation $x^n = Dy^2 + 1$"}, Acta Arith. {\bf 113} (2004), no. 1, 69--76.
\bibitem{HE99} C. Heuberger and M. Le, {\it On the generalized Ramanujan-Nagell equation $x^2 + D = p^Z$}, J. Number Theory {\bf 78} (1999), no. 2, 312--331.

\bibitem{HC19} A. Hoque and K. Chakraborty, {\it Divisibility of class numbers of certain families of quadratic fields}, J. Ramanujan Math. Soc. {\bf 34} (2019), no. 3, 281--289. 
\bibitem{HU87} L. K. Hua, {\it Introduction to number theory}, Springer-Verlag, New York, 1982.
\bibitem{IZ18} Y. Iizuka, {\it On the class number divisibility of pairs of imaginary quadratic fields}, J. Number Theory {\bf 184} (2018), 122--127.
\bibitem{IK21} Y. Iizuka, Y. Konomi and S. Nakano, {\it An application of the arithmetic of elliptic curves to the class number problem for quadratic fields}, Tokyo J. Math. (2021). doi: \url{10.3836/tjm/1502179314}
\bibitem{KO02} T. Komatsu, {\it An infinite family of pairs of quadratic fields $\mathbb{Q}(\sqrt{D})$ and $\mathbb{Q}(\sqrt{mD})$ whose class numbers are both divisible by $3$}, Acta Arith. {\bf 104} (2002), 129--136.

\bibitem{KO17} T. Komatsu, {\it An infinite family of pairs of imaginary quadratic fields with ideal classes of a given order}, Int. J. Number Theory {\bf 13} (2017), no. 2, 253--260.

\bibitem{KP21} S. Krishnamoorthy and S. Pasupulati, {\it Note on the $p$-divisibility of class numbers of an infinite family of imaginary quadratic fields}, Glasgow Math. J. (2021). doi:    \url{10.1017/S001708952100015X}
\bibitem{LO09} S. R. Louboutin, {\it On the divisibility of the class number of imaginary quadratic number fields}, Proc. Amer. Math. Soc. {\bf 137} (2009), no. 12, 4025--4028.

\bibitem{RM98}M. R. Murty, {\it The ABC conjecture and exponents of class groups of quadratic fields}, Contemp. Math. {\bf 210} (1998), 85--95.

\bibitem{NA22} T. Nagell, {\it \"{U}ber die Klassenzahl imagin\"{a}r quadratischer, Z\"{a}hlk\"{o}rper}, Abh. Math. Sem. Univ. Hambg. {\bf 1} (1922), 140--150.

\bibitem{SC32} A. Scholz, {\it \"Uber die Beziehung der Klassenzahlen quadratischer K\"orper zueinander}, J. Reine Angew. Math. {\bf 166} (1932), 201--203.

\bibitem{SI29} C. L. Siegel, {\it Uber einige Anwendungen Diophantischer Approximationen}, Abh. Preuss. Akad. Wiss. Phys. Math. Kl. {\bf 1} (1929), 1-70; Ges. Abh., Band {\bf 1}, 209--266.
\bibitem{SO00} K. Soundararajan, {\it Divisibility of class numbers of imaginary quadratic fields}, J. London Math. Soc. {\bf 61} (2000), 681--690.

\bibitem{VO95} P. M. Voutier, {\it Primitive divisors of Lucas and Lehmer sequences}, Math. Comp. {\bf 64} (1995), 869--888.
\bibitem{XC20} C. -F. Xie and C. F. Chao, {\it On the divisibility of class numbers of imaginary quadratic fields $\left(\mathbb{Q}(\sqrt{D}),\mathbb{Q}(\sqrt{D+m})\right)$}, Ramanujan J. {\bf 53} (2020), 517--528.
\bibitem{YA70} Y. Yamamoto, {\it On unramified Galois extensions of quadratic number fields}, Osaka J. Math. {\bf 7} (1970), 57--76.

\end{thebibliography}
\end{document}